\documentclass[letterpaper]{article}

\usepackage[margin=2.54cm]{geometry}
\usepackage{setspace}

\usepackage{amsthm}
\usepackage{amsmath}
\usepackage{amsfonts}
\usepackage{amssymb}
\usepackage{amscd}
\usepackage{mathrsfs}
\usepackage{bm}
\usepackage{url,hyperref}

\theoremstyle{plain}

\theoremstyle{definition}

\newcommand{\Q}{\mathbb{Q}}

\newcommand{\pa}{{\rm pa}}

\def\ci{\perp\!\!\!\perp}

\title{Graphical models in Macaulay2}
\author{Luis David  Garc\'{\i}a-Puente\footnote{Department of Mathematics and Statistics,         Sam Houston State University,          Huntsville, TX\ 77341, \href{mailto:lgarcia@shsu.edu}{lgarcia@shsu.edu}},
 Sonja Petrovi\'c\footnote{Department of Statistics,       Pennsylvania State University,           University Park, PA\ 16802, \href{mailto:petrovic@psu.edu}{sonja@psu.edu}}, 
 Seth Sullivant\footnote{Department of Mathematics,   North Carolina State University, Raleigh, NC\ 27695, \href{mailto:smsulli2@ncsu.edu}{smsulli2@ncsu.edu}}}

\date{\today}

\begin{document}

\maketitle

\begin{abstract}
The \emph{Macaulay2} package {\tt GraphicalModels} contains algorithms 
for the algebraic study of graphical models associated to undirected, 
directed and mixed graphs, 
and associated collections of 
conditional independence statements. 
Among the algorithms implemented are procedures for computing
the vanishing ideal of graphical models, for generating
conditional independence ideals of families of independence statements
associated to graphs, and for checking for identifiable parameters
in Gaussian mixed graph models.
These procedures can be used to study fundamental 
problems about graphical models. 
\end{abstract}

\section{Graphical models}\label{intro}

A graphical model is a statistical model associated to a graph,
where the nodes of the graph represent random variables and the
edges of the graph encode relationships between the random variables.
Graphical models are an important class of statistical models used
in many applications (see standard textbooks \cite{Lauritzen, Whitaker})
because of their ability to model complex interactions between
several random variables, by specifying
interactions using only local information about connectivity
between the vertices in a graph.

There are two natural ways to specify a graphical model, through 
either conditional independence statements specified by the graph
or via a parametric representation (often called a ``factorization'').
Every distribution that factors according to the graph
satisfies the conditional independence statements implied by the 
graph. This leads to the question:  Which distributions satisfy the conditional
independence statements implied by the graph, but do not factor?

Once we specify the types of random variables under consideration
(e.g., discrete random variables or Gaussian random variables) it
is possible to address the questions in the preceding paragraph
using (computational) algebraic geometry.  Indeed, in these cases,
the set of all probability distributions satisfying a family of 
conditional independence constraints is a semialgebraic set.
For discrete random variables, that semialgebraic set is a subset
of the probability simplex, and can be represented by a certain
homogeneous ideal generated by quadrics.  For Gaussian random variables,
this set of distributions corresponds to a semialgebraic subset
of the cone of positive definite matrices.  Similarly,
the parametrized family of probability distributions also
is a semialgebraic set (of the probability simplex for discrete
random variables, and of the cone of positive definite matrices
for Gaussian random variables).  This algebraic perspective
has been studied by different authors  
\cite{ GSS, GMS, S}, and the book \cite{DSS}
provides  details.

The \emph{Macaulay2} package {\tt GraphicalModels} allows the user to compute
the ideals of conditional independence statements for any collection
of statements for discrete or Gaussian
random variables.  It can also compute the vanishing ideal of 
a graphical model in these cases.  A number of auxiliary functions
are useful for doing further analyses of graphical models.

For example, consider the directed acyclic graph $G$ with  five
vertices $\{a,b,c,d,e\}$ and edge
set $\{a \to d, b \to d, c \to d, c \to e, d \to e\}$.
The following commands compute the associated conditional
independence ideal for the set of global Markov statements,
$CI_{{\rm global}(G)}$, and 
the vanishing ideal $I_{G}$ of the Gaussian graphical model on $G$.
In the following examples,  version 1.4 of \emph{Macaulay2} was used. 


\begin{verbatim}
i1 : loadPackage "GraphicalModels";
i2 : G = digraph{{a,d},{b,d},{c,{d,e}},{d,e}}; 
i3 : R = gaussianRing G
o3 : PolynomialRing
i4 : gens R
o4 = {s   , s   , s   , s   , s   , s   , s   , s   , s   , s   , s   , s   , 
       a,a   a,b   a,c   a,d   a,e   b,b   b,c   b,d   b,e   c,c   c,d   c,e   
      s   , s   , s   } 
       d,d   d,e   e,e  
i5 : I = conditionalIndependenceIdeal(R,globalMarkov(G));
i6 : J = gaussianVanishingIdeal(R);
i7 : flatten degrees J
o7 = {1, 1, 1, 2, 3, 3}
\end{verbatim}

 {\tt GraphicalModels} uses the package {\tt Graphs} and a number of fundamental constructs and relationships 
associated with graphs.  First we create a polynomial ring that contains the entries of
the covariance matrix $\Sigma$ of a jointly normal random vector
as its indeterminates.  Information about the underlying graph is stored in the polynomial ring. 
Hence some methods take just a ring as input, but require that it be created with {\tt gaussianRing}, or {\tt markovRing} in the discrete case.

For directed acyclic graphs it is known that 
$V(CI_{{\rm global}(G)}) \cap PD_{m}  =  
V(I_{G}) \cap PD_{m}$,
in particular, the set of positive definite matrices
satisfying the conditional independence constraints equals the set of covariance matrices in the image of the parametrization.
Unfortunately, this does not imply that $CI_{{\rm global}(G)} 
=I_{G}$.  In the case of Gaussian random variables, a larger ideal,
the trek ideal $T_{G}$, generated by all subdeterminants of the 
covariance matrix that vanish on the model, and satisfying
$CI_{{\rm global}(G)} \subseteq T_{G} \subseteq
I_{G}$ is sometimes equal to $I_{G}$ (see \cite{STD}), as the following example shows. 

\begin{verbatim}
i8 : isSubset(I,J)
o8 = true
i9 : I == J
o9 = false
i10 : J == trekIdeal(R,G)
o10 = true
\end{verbatim}

Similar computations can also be performed for graphical models
with discrete random variables,
and with other graph families.  The mathematical explanation of
these graphical models and their associated ideals
appear in the remaining sections.

\section{Computing conditional independence ideals}\label{CI}

Conditional independence constraints on discrete or Gaussian random variables
 translate to rank conditions on certain matrices
associated to the probability densities.  We briefly explain
these constructions here and how to generate these constraints
in \emph{Macaulay2} using {\tt GraphicalModels}. See \cite[Ch.~3]{DSS}
for more detail.

Let $X = (X_{1}, \ldots, X_{n})$ be a discrete random
vector where each random variable $X_{i}$ has state space
$[d_{i}] = \{1,2, \ldots, d_{i} \}$. Let $d = (d_{1}, \ldots, d_{n})$.  
A probability distribution
for $X$ is a tensor in $\mathbb{R}^{d_{1}}\otimes \cdots \otimes
\mathbb{R}^{d_{n}}$, all of whose 
coordinates are nonnegative and sum to one.
The set of all such distributions is the probability simplex $\Delta_{d}$.
Let $p_{i_{1}\cdots i_{n}} = {\rm P}(X_{1} = i_{1}, \ldots, X_{n} = i_{n})$
denote the probability of a primitive event.  The polynomial
ring in these quantities is created using the command ${\tt markovRing}$.

\begin{verbatim}
i11 : d = (2,3,2); R = markovRing d
o11 : PolynomialRing
i12 : gens R
o12 = {p     , p     , p     , p     , p     , p     , p     , p     , p     ,
        1,1,1   1,1,2   1,2,1   1,2,2   1,3,1   1,3,2   2,1,1   2,1,2   2,2,1 
      p     , p     , p     }
       2,2,2   2,3,1   2,3,2
\end{verbatim}

For $A \subseteq [n]$, let $X_{A} =(X_a)_{a\in A}$ be the subvector
indexed by $A$.  Let $A,B,C$ be disjoint subsets of $[n]$.  
The conditional independence statement
$X_{A} \ci X_{B} | X_{C}$ holds if and only if the conditional
distribution satisfies
$$
{\rm P}(X_{A} = i_{A}, X_{B} = i_{B} | X_{C} = i_{C}) =
{\rm P}(X_{A} = i_{A} | X_{C} = i_{C})\cdot {\rm P}( X_{B} = i_{B} | X_{C} = i_{C})$$
for all $i_{A}, i_{B}, i_{C}$.   This translates into vanishing
$2\!\times\! 2$ minors of certain matrices in the probabilities 
$p_{i_{1}\cdots i_{n}}$.  Those matrices are computed with the 
function {\tt markovMatrices},
and the ideal generated by the $2 \!\times\! 2$ minors is computed with 
{\tt conditionalIndependenceIdeal}.  In the following example,
the two conditional independence statements are $X_{1} \ci X_{2} | X_{3}$
and $X_{1} \ci X_{3} \, \,  (:=  X_{1} \ci X_{3} | X_{\emptyset})$. 
The ideal of vanishing minors has $7$ quadratic generators. 

\begin{verbatim}
i13 : s = {{{1},{2},{3}}, {{1},{3},{}}}
i14 : compactMatrixForm=false;
i15 : markovMatrices(R,s)
o15 = {| p       p       p      |, | p       p       p      |, 
       |  1,1,1   1,2,1   1,3,1 |  |  1,1,2   1,2,2   1,3,2 | 
       |                        |  |                        |  
       | p       p       p      |  | p       p       p      |  
       |  2,1,1   2,2,1   2,3,1 |  |  2,1,2   2,2,2   2,3,2 |
             
       | p      + p      + p       p      + p      + p      |} 
       |  1,1,1    1,2,1    1,3,1   1,1,2    1,2,2    1,3,2 |
       |                                                    |
       | p      + p      + p       p      + p      + p      |
       |  2,1,1    2,2,1    2,3,1   2,1,2    2,2,2    2,3,2 |  
i16 : I = conditionalIndependenceIdeal(R,s);
o16 : Ideal of R
i17 : flatten degrees  I 
o17 = {2, 2, 2, 2, 2, 2, 2}
\end{verbatim}

In  statistics literature, there are three main lists of conditional independence statements 
associated to a graph $G$ whose nodes correspond to random variables. 
 For example, the list of \emph{local Markov statements} of an undirected graph $G$
is the set of conditional independence statements of the form
$X_{i} \ci X_{V \setminus \left\{i \cup N(i)\right\}} | X_{N(i)}$, where $N(i)$ is the set of neighbors of $i$ in the graph $G$.
The methods {\tt pairMarkov}, {\tt localMarkov}, and {\tt globalMarkov}
compute the pairwise, local, and global Markov statements, respectively, 
for both directed and undirected graphs.

\begin{verbatim}
i18 : G = graph{{1,2},{2,3},{3,4},{4,5},{1,5}};
i19 : localMarkov G
o19 = {{{1}, {3, 4}, {5, 2}}, {{1, 2}, {4}, {5, 3}}, {{1, 5}, {3}, {4, 2}}, 
      {{2, 3}, {5}, {4, 1}}, {{2}, {4, 5}, {1, 3}}}
\end{verbatim}

For example, the first conditional independence statement produced
is $X_{1} \ci (X_{3}, X_{4}) | (X_{2}, X_{5})$.  In the 
context of conditional independence, the graphical model
consists of all distributions satisfying one of these collections
of independence statements associated to the graph $G$.

A Gaussian random vector, $X = (X_{1}, \ldots, X_{n}) 
\sim \mathcal{N}(\mu, \Sigma)$, is an
$n$-dimensional random vector with state space $\mathbb{R}^{n}$ and
density function
$$
f(x)  =  \frac{1}{(2 \pi)^{n/2}(\det \Sigma)^{1/2} } 
\exp\left( - \frac{1}{2} (x- \mu)^{T} \Sigma^{{-1}} (x - \mu) \right), 
$$
where $\mu \in \mathbb{R}^{n}$ and $\Sigma = (\sigma_{s,t}) \in PD_{n}$, the cone of 
$n \!\times\! n$ symmetric positive definite matrices.  The Gaussian random vector
$X$ satisfies the conditional independence statement $X_{A} \ci X_{B} | X_{C}$
if and only if the submatrix 
$\Sigma_{A \cup C, B \cup C}  := 
(\sigma_{s,t})_{s \in A \cup C, t \in B \cup C}$ has rank $\leq \#C$.
Hence the set of all Gaussian random vectors satisfying a given
collection of conditional independence statements yields a 
subset of $PD_{n}$ that can be studied via
a determinantal conditional independence ideal in the polynomial ring
$\Q[\sigma_{s,t}: s,t \in [n]]$.  This ring
is generated using the command {\tt gaussianRing}.
Computations involving conditional independence ideals with
Gaussian random variables were exemplified in Section \ref{intro}.

\section{Computing the vanishing ideal of a model}

The fact that graphical models can be described in two possible ways (either by a recursive
factorization of probability distributions or by conditional independence
statements) corresponds to the algebraic principle that varieties can be presented either
parametrically or implicitly. The \emph{vanishing ideal} of a model is the set
of homogeneous polynomial relations in the joint probability distributions (for discrete
random variables) or in the variance-covariance parameters (for Gaussian
random variables).
{\tt GraphicalModels} has the capability of computing the vanishing ideals
of graphical models on directed graphs (for discrete random variables) and also
of graphical models on directed, undirected, or  mixed
graphs (for Gaussian random variables).  The vanishing ideal
of an undirected graphical model for discrete random variables is
a toric ideal and
should be computed using the \emph{Macaulay2} package {\tt FourTiTwo}.

The method
{\tt discreteVanishingIdeal} implements this capability for graphical
models on discrete random variables. For a directed acyclic graph $G$ on discrete random variables, the graphical model consists
of all distributions satisfying the recursive factorization
property 
\[p(X = i) = \prod_{v} {\rm P}(X_{v} = i_{v} | X_{\pa(v)} = i_{\pa(v)}),\]
where the product runs over all vertices $v$ of $G$ and $\pa(v)$ is the set of
parents of $v$. 
Our implementation of this method does not compute the kernel of the
corresponding ring map. Instead, the vanishing ideal is
computed recursively using the factorization 
\[{\rm P}(X = i) = {\rm P}(X_{1} = i_{1},\dots, X_{n-1} = i_{n-1})\cdot {\rm P}(X_{n} = i_{n} | X_{\pa(n)} = i_{\pa(n)}),\]
where $1, \ldots, n$ is a topological ordering of
the vertices of the directed acyclic graph $G$.

The following example computes the vanishing ideal of the graphical model $1\to
2 \to 3\to 4$ on four binary random variables. The vanishing ideal is
minimally generated by 20 quadratic binomials.

\begin{verbatim}
i20 : G = digraph {{1,{2}}, {2,{3}},{3,{4}},{4,{}}};
i21 : R = markovRing (2,2,2,2);
i22 : I = discreteVanishingIdeal (R,G);
o22 : Ideal of R
i23 : betti mingens I
             0  1
o23 = total: 1 20
          0: 1  .
          1: . 20
\end{verbatim}

According to \cite{GSS}, the vanishing ideal of a graphical model on discrete
random variables is the \emph{distinguished
  component} of the conditional independence ideal described by the Markov
statements of the model. For the directed path in our previous example, the
conditional independence ideal of the local Markov statements 
is a radical ideal with 3 associated primes. However, since $G$ is a directed
tree, the conditional independence ideal of the global Markov statements is a
prime ideal and it equals the vanishing ideal of $G$.

 \begin{verbatim}
i24 : J = conditionalIndependenceIdeal (R, localMarkov G);
o24 : Ideal of R
i25 : I == J
o25 = false
i26 : K = conditionalIndependenceIdeal (R, globalMarkov G);
o26 : Ideal of R
i27 : I == K
o27 = true
\end{verbatim}

The method {\tt gaussianVanishingIdeal} computes the vanishing ideal of a
Gaussian graphical model on a graph, digraph, or mixed graph. If $G$ is a mixed
graph, 
$I_{G}$ is the vanishing ideal of the parametrization 
\[\Sigma = (I - \Lambda)^{-T}\Psi(I-\Lambda)^{-1},\] 
where $\Sigma$ is the variance-covariance matrix, $\Lambda$ is the strictly
upper triangular matrix with $\Lambda_{ij} = \lambda_{ij}$ if $i\to j$ is a
directed edge in $G$ and 0 otherwise, and $\Psi$ is a symmetric positive definite
matrix of parameters $\psi_{ij}$ with zeros in each entry $\Psi_{ij}$ if there is no
bidirected edge in $G$ between $i$ and $j$, and $i\neq j$.

The following example computes the vanishing ideal of the Gaussian graphical
model  on the mixed graph with directed edges
$\{1\to 2, 1\to 3, 2\to 3, 3\to 4\}$ and bidirected edges $\{1 \leftrightarrow 2, 2 \leftrightarrow
4\}$. This ideal is a principal ideal generated by one quartic polynomial with 8
terms. This ideal is not \emph{determinantal}, i.e., it is not generated by the
determinantal equations defining 
the trek ideal, which in this case is the zero ideal.

\begin{verbatim}
i28 : G = mixedGraph(digraph {{1,{2,3}},{2,{3}},{3,{4}}},bigraph {{1,2},{2,4}});
i29 : R = gaussianRing G;
i30 : I = gaussianVanishingIdeal R;
o30 : Ideal of R
i31 :  flatten degrees  I 
031 = {4}
i32 : J = trekIdeal (R,G)
o32 = 0
\end{verbatim}

An important problem in these models consists in finding which parameters are
\emph{identifiable} or \emph{generically identifiable}, see \cite{GPSS}. The method
\texttt{identifyParameters} can be used to solve the identifiability problem for
Gaussian graphical  models on mixed graphs (also known as \emph{structural
  equation models}). The following example shows that the parameter
$\psi_{24}$ is generically identifiable by the formula $\psi_{24} =
(\sigma_{13}\sigma_{24} - \sigma_{14}\sigma_{23})/ \sigma_{13}$. 

\begin{verbatim}
i33 : H = identifyParameters R;
i34 : H#(p_(2,4))_0
o34 = p   s    + s   s    - s   s
       2,4 1,3    1,4 2,3    1,3 2,4
\end{verbatim}

In this model there are three non-generically identified parameters.  \texttt{identifyParameters} produces a hash table whose entries
are indexed by the parameters and contain ideals that can be used to find
explicit rational functions for every parameter that is generically
identifiable.

\section*{Acknowledgments }
The following people generously contributed their time to the
development of the package:  Alexander Diaz, Shaowei Lin, David Murrugarra, and Mike 
Stillman.  Work on the package was carried out during the 2010 and 2011
\emph{Macaulay2} workshops, which were partially supported by the US National 
Science Foundation and the Institute for Mathematics and Its
Applications. The authors also thank the anonymous referees for their helpful comments and 
suggestions that improved not only this article but also the GraphicalModels package.

LGP was partially supported by a 2012 SHSU Faculty Research Grant (29001).
SP was partially supported by Grant \#FA9550-12-1-0392 from the U.S. Air Force Office of Scientific Research (AFOSR) and the Defense Advanced Research Projects Agency (DARPA).
SS was partially supported by the US National Science Foundation (DMS 0954865) and the David and Lucille Packard Foundation.

\begin{spacing}{0.5}

\end{spacing}

\end{document}